\documentclass[psamsfonts, 12pt, a4, leqno]{amsart}

\usepackage{amssymb, latexsym}

\newtheorem{Theorem}{Theorem}
\newtheorem{Corollary}[Theorem]{Corollary}
\newtheorem{Lemma}[Theorem]{Lemma}
\newtheorem{Proposition}[Theorem]{Proposition}

\newcommand{\A}{\operatorname{area}(X)}

\usepackage
{hyperref}

\begin{document}

\title[On primeness of the Selberg zeta-function] {On primeness of the Selberg zeta-function}

\author{ Ram\={u}nas Garunk\v{s}tis}
\address{Ram\={u}nas Garunk\v{s}tis \\
Institute of Mathematics \\ Faculty of Mathematics and Informatics\\ Vilnius University \\
Naugarduko 24, 03225 Vilnius, Lithuania}
\email{ramunas.garunkstis@mif.vu.lt}
\urladdr{www.mif.vu.lt/~garunkstis}

\author{ J\"orn Steuding}
\address{J\"orn Steuding \\
Department of
Mathematics, W\"urzburg University \\
Am Hubland,
97\,218 W\"urzburg, Germany}
\email{steuding@mathematik.uni-wuerzburg.de}
\urladdr{}

\setcounter{equation}{0}

\date{}
\maketitle

\begin{abstract}
In this note we prove that the Selberg zeta-function associated to a compact Riemann surface is pseudo-prime and right-prime in the sense of a decomposition. 
\end{abstract}


\section{Introduction}

Let $s=\sigma+it$ be a complex variable and $X$ a compact Riemann
surface of genus $g \geq 2$. The surface $X$ can be regarded as a
quotient $\Gamma \backslash H$, where $\Gamma \subset\text{PSL}(2,
\mathbb{R})$ is a strictly hyperbolic Fuchsian group and $H$ is the
upper half-plane of $\mathbb{C}$.  Then the Selberg zeta-function
associated with $X=\Gamma \backslash H$ is defined by (see
Hejhal~\cite[\S 2.4, Definition 4.1]{hejhal76})
\begin{equation}\label{eq:1}
  Z(s) = \prod_{\{P_0 \}}\prod_{k = 0}^\infty(1 - N(P_0)^{-s-k}).
\end{equation}
Here $\{ P_0 \}$ is the conjugacy class of a primitive hyperbolic element $P_0$ of $\Gamma$ and $N(P_0)
=\alpha^2$ if the eigenvalues of $P_0$ are $\alpha$ and $\alpha^{-1}$
with $|\alpha| > 1$. Equation \eqref{eq:1} defines the Selberg
zeta-function in the half-plane $\sigma> 1$. The function $Z(s)$ can
be extended to an entire function of order~2 (see \cite[\S 2.4,
Theorem 4.25]{hejhal76}). 

The Selberg zeta-function $Z(s)$ has so-called trivial zeros at integers $s=-n$, $n\ge1$, of multiplicity $(2g-2)(2n+1)$; at $s=0$ with multiplicity $2g-1$; and at $s=1$ with multiplicity $1$. There are further, so-called nontrivial zeros on the critical line $\sigma=1/2$ with
at most finitely many exceptions of zeros on the real segment $0<s<1$
(see \cite[\S 2.4, Theorem 4.11]{hejhal76} and Randol
\cite{randol74}). The nontrivial zeros $s_j=1/2\pm it_j$ correspond to
eigenvalues
\begin{align}\label{eigenvalueszeros}
\lambda_j=s_j(1-s_j)=1/4+t_j^2>0
\end{align}
of the hyperbolic Laplacian $\Delta$ on $X=\Gamma \backslash H$
(see \cite[\S 2.4, Theorem 4.11]{hejhal76}. The imaginary parts of nontrivial zeros (and, more generally, imaginary parts of $a$-points for any given complex number $a$) of $Z(s)$ are uniformly distributed modulo one (see \cite{gss14}). 

Moreover, the Selberg
zeta-function satisfies the following functional equation (\cite[\S 2.4, Theorem 4.12]{hejhal76})
\begin{equation}\label{eq:2}
  Z(s) = X(s) Z(1 - s),
\end{equation}
where
\begin{equation*}
  X(s) = \exp \left( 4 \pi (g - 1) \int_0^{s - 1/2} v \tan (\pi v) \, \text{d}v \right).
\end{equation*}

We consider the compositions $F(z) = f(h(z))$, where $F$, $f$, and $h$ are meromorphic functions on $\mathbb C$. We always can achieve that $f$ is meromorphic and $h$ is an entire function. However, if $F$ has infinitely many poles and $f$ is a rational function, then $h$ is meromorphic. We arrive at the following definition.
\medskip

\noindent {\bf Definition} (Gross \cite{gro} or Chuang and Yang \cite[Section 3.2]{cc}).  Let $F$ be a meromorphic function. Then an expression
\begin{align}\label{factorization}
F(z) = f(h(z))
\end{align} 
where $f$ is meromorphic and $h$ is entire ($h$ may be meromorphic when $f$ is a rational function) is called a {\it decomposition} of $F$ with $f$ and $h$ as its left and right  components 
respectively. $F$ is said to be  {\it prime} in the sense of a decomposition if for every representation
of $F$ of the form \eqref{factorization} we have that either $f$ or $h$ is linear. If every representation
of $F$ of the form \eqref{factorization} implies that $f$ is rational or $h$ is a polynomial ($f$ is linear
whenever $g$ is transcendental, $g$ is linear whenever $f$ is transcendental), we say
that $F$ is {\it pseudo-prime} ({\it left-prime}, {\it right-prime}) in the sense of a decomposition. Note that here the terminology is slightly changed. In \cite{gro} and \cite{cc} the notions {\it factorization} and  {\it factor} instead of the corresponding notions {\it decomposition} and {\it component} were used. Further we use the shorter wording `function is prime (pseudo-prime)' instead of `function is prime (pseudo-prime) in the sense of a decomposition'.
\medskip

The first example of a prime function in the literature is $F(z)=e^z+z$ (see Rosenbloom \cite{ros} and Gross \cite{gro}). Liao and Yang \cite{ly} proved that both, the Gamma function and the Riemann zeta-function are prime functions. Here we are concerned with another type of zeta-function, namely the Selberg zeta-function. We expect that it is a prime function. Our attempts led to the following

\begin{Theorem}\label{prime}
The Selberg zeta-function $Z$ associated with a compact Riemann surface of genus $g$ is pseudo-prime and right-prime. 

Moreover, if $Z(s) = f(h(s))$, where $f$ is rational and $h$ is meromorphic, then $f$ is a polynomial of degree $k$, where $k$ divides $2g-2$, and $h$ is an entire function.
\end{Theorem}

\noindent Indeed, the polynomial $f$ in this theorem cannot be of the form $f(z)=cz^k$, $k\ge2$, since $Z(s)$ has a simple zero at $s=1$. In the proof of Theorem \ref{prime} the difference of the growth of $Z(s)$ in right and left half-planes of $\mathbb C$ is used (a topic for which we refer to \cite{GarunkstisGrigutis16}).

The following consequence of our theorem might be of independent interest.

\begin{Corollary}\label{z1z2strong} 
Let $Z_1$  and $Z_2$ be the Selberg zeta-functions associated with compact Riemann surfaces $X_1$ and $X_2$, respectively. For $j=1,2$, denote by $S_j$ the set of nontrivial zeros of $Z_j(s)$. Assume that $h$ is an entire function satisfying $h^{-1}(S_2)= S_1$. Then $h$ is the identity and $S_1=S_2$.
\end{Corollary}

Theorem \ref{prime} and Corollary \ref{z1z2strong} are proved in the next section. In Section \ref{cr} we discuss a hypothetical decomposition of $Z$ with a quadratic polynomial and an entire function.

\section{Proofs}

We begin with a 
\medskip

\noindent {\bf Definition.}  Let $E\subset \mathbb C$. If $\theta\in[0,2\pi]$ is an accumulation point of $S=\{\arg z : z\in E\}$, then we call $\{z\,:\,\arg z=\theta\}$ an {\it accumulation line} of $E$.
\medskip

Liao and Yang \cite{ly} used the following lemma to prove the pseudo-primeness of $\zeta$. From this lemma we shall derive that $Z$ is also pseudo-prime. 

\begin{Lemma}\label{accumlines}
Let $a_1,a_2$ be arbitrary distinct complex numbers or $\infty$ and let $f$ be a meromorphic function of finite order. Assume that the number of the accumulation lines of $E=\{z : f(z)=a_j\in \overline{\mathbb C}, j=1,2\}$ is finite. Then $f$ is pseudo-prime.
\end{Lemma}

\proof Lemma \ref{accumlines} is proved in \cite[p. 141]{chu}. \endproof

The following lemma describes the asymptotic behavior of the factor $X(s)$ in the functional equation (\ref{eq:2}).
\begin{Lemma}\label{ft}
  For $t\ge1$,
  \begin{align*}
    X(s) = &\exp \left( 2 \pi i (g-1)\left(s - \frac{1}{2} \right)^2 +\frac{\pi i (g-1)}{6} \right.\\
    \nonumber &\quad\qquad \left.+O \left(\frac{t}{\text{e}^{2 \pi t}}+ \frac{(\sigma - 1/2)^2}{\text{e}^{2 \pi t}} + \frac{(\sigma -1/2)t}{\text{e}^{2 \pi t}}\right)
    \right)\qquad (t\to\infty)
  \end{align*}
  uniformly in $\sigma$.
\end{Lemma}

\proof This is Lemma 1 in \cite{garunkstis12}. \endproof

\begin{Proposition}\label{propo}
The Selberg zeta-function associated to a compact Riemann surface is pseudo-prime and right-prime.
\end{Proposition}

\proof By Lemma \ref{accumlines} (applied with $a_1=0$ and $a_2=\infty$) the Selberg zeta-function $Z(s)$ is pseudo-prime.

Next we show that $Z$ is a right-prime function. Assume that 
\begin{align*}
Z(s)=f(h(s))
\end{align*}
where  $f$ is entire and $h$ is a polynomial. We need to prove that then $h$ is a linear function. For this we shall use growth properties of the Selberg zeta-function. Let
$
N(P_{00})=\min_{P_0}\{N(P_0)\}.
$
By definition of the norm we have $N(P_{00})>1$. In \cite[Proposition 4.13]{hejhal76} it is proved that
\begin{align}\label{1}
Z(s)=1+O({N(P_{00})}^{-\sigma})\quad(\sigma\to\infty),
\end{align}
uniformly in $\sigma\ge2$. Combining this with the functional equation \eqref{eq:2} and Lemma \ref{ft}, we observe that, for any $\varepsilon>0$,
\begin{align}\label{onleft}
\lim_{|s|\to\infty} Z(s)=\infty
\end{align}
uniformly in 
$$A:=\{s : \pi/2+\varepsilon\le\arg s\le\pi-\varepsilon\quad\text{or}\quad \pi+\varepsilon\le\arg s\le3\pi/2-\varepsilon\}.$$
Let $h(s)=a_ds^d+a_{d-1}s^{d-1}+\ldots+a_1s+a_0$ with $a_d\neq 0$. For the polynomial $h$ we consider the preimages $\ell_1$, \dots, $\ell_d$ of the half-line 
$$\ell=\left\{s : \arg s=\frac{\pi}{2}-\arg a_d\right\}.$$
We number the preimages such that near to infinity the curve $\ell_j$ is close to the half-line 
$$L_j:=\left\{s : \arg s=\frac{\pi}{2d}-\frac{\arg a_d}{d}+\frac{2j\pi}{d}\right\},$$
where $j\in\{1,\dots, d\}$. Note that $\arg L_j\ne \pi$, $j\in\{1,\dots, d\}$.  Thus, if $d\ge2$, then there are indices $p,q\in\{1,\dots, d\}$ such that $L_p$ lies in the half-plane $\sigma>2$ (except for a finite part of $L_p$) and $L_q$ lies in the set $A$. Therefore, by formulas \eqref{1} and \eqref{onleft}, we have
\begin{align}\label{imposible}
\lim_{{|s|\to\infty\atop s\in\ell_p}}Z(s)=1\quad\text{and}\quad\lim_{{|s|\to\infty\atop s\in\ell_q}}Z(s)=\infty.
\end{align}
>From other side, by $Z(s)=f(h(s))$ and $h(\ell_p)=\ell=h(\ell_q)$ we see that $Z(\ell_p)=Z(\ell_q)$. The last equality contradicts equation \eqref{imposible}. Thus $d=1$. This proves Proposition \ref{propo}.
\endproof

\begin{Lemma}\label{pseaudoprimepolynomial}
If $F$ is entire right-prime, $F(z)=f(h(z))$ with rational $f$ and entire $h$, then $f$ is a polynomial or $h(z)=\exp(az+b)-w$ for some complex numbers $a\ne0$, $b$, and $w$.
\end{Lemma}

\proof By Picard's theorem $f$ has at most one pole. Suppose $f(z)$ has one pole at $z=w$. Then $h(z)$ omits the value $w$. Hence, there is an entire function $p$, such that $h(z)=\exp(p(z))-w$. Then $p$ is a linear function since $F$ is a right-prime function. \endproof

Let
\begin{align}\label{Fdef}
F(s)=\A \int_0^sz\tan(\pi z) dz,
\end{align}
where the integration is along the straight line segment joining the origin to $s$ if $s$ is not on the real line; otherwise, when $s$ is on the real line, and not one of the
points $\pm1/2$, $\pm3/2$, $\pm5/2$,\dots, we define $F(s)$ by the requirement of
continuity as $s$ is approached from the upper half-plane (compare to the definition of the function $\Phi(s)$ in Randol \cite[proof of Lemma 2]{randol78}).

We shall find an antiderivative of $z\tan(\pi z)$ in \eqref{Fdef}. The dilogarithm function is the function defined by the power series
\begin{align}\label{powerseriesli}
\operatorname{Li}_2(s)=\sum_{n=1}^\infty\frac{s^n}{n^2}\qquad\text{for}\ |s|\le1. 
\end{align}
Analytic continuation of the dilogarithm is given by
\begin{align}\label{continuationli}
\operatorname{Li}_2(s)=-\int_0^s\log(1-z)\frac{dz}{z} \qquad\text{for}\ s\in\mathbb C\setminus[1,\infty).
\end{align}
For $t>0$, let 
$$P(s)=\A \left(\frac{is^2}{2} -\frac{s}\pi\log\left(1+e^{2i\pi s}\right)+\frac{i}{2\pi^2}\operatorname{Li}_2(-e^{2i\pi s})+\frac{i}{24}\right),$$
where the principal branch of the logarithm is chosen.
In view of $\operatorname{Li}_2(-1)=-\pi^2/12$ (see formulas (1.8) and (1.9) in Lewin \cite{lew}) it follows that $\lim_{s\to0} P(s)=0=F(0)$. Taking into account the expressions \eqref{Fdef} and \eqref{continuationli} we   obtain that $P'=F'$. Thus, for $t>0$,
\begin{align}\label{ft}
F(s) = \A \left(\frac{is^2}{2} -\frac{s}\pi\log\left(1+e^{2i\pi s}\right)+\frac{i}{2\pi^2}\operatorname{Li}_2(-e^{2i\pi s})+\frac{i}{24}\right).
\end{align}

\begin{Lemma}\label{sn12}
For $s$ satisfying $|s-n|=1/2$ with a negative integer $n$, we have $Z(s)\to\infty$ as $n\to-\infty$.
\end{Lemma}
 
\proof In view of the functional equation \eqref{eq:2} and formula \eqref{1} we need to show that 
\begin{align}\label{xtoinfty}
X(s)\to\infty 
\end{align} 
for $|s-n|=1/2$, as $n\to-\infty$.  We have $|X(s)|=\exp(\Re F(s-1/2))$. Then \eqref{ft} together with  expressions \eqref{powerseriesli},
\begin{align*}
\Re \left(i(s-1/2)^2/2\right)=-t(\sigma-1/2),
\end{align*}
the fact that 
\begin{align*}
\Re\left(\frac{s-1/2}{\pi}\log(1+e^{2i\pi(s-1/2)})\right)=\frac{\sigma-1/2}{\pi}\log|1-e^{2i\pi s}|+O(1)\quad (\sigma\to-\infty),
\end{align*}
uniformly in $0<t\le1$, yield
\begin{align*}
|X(s)|=\exp\left(\A(1/2-\sigma)(t+\frac{1}{\pi}\log|1-e^{2i\pi s}|)+O(1)\right)\quad (\sigma\to-\infty),
\end{align*}
uniformly in $0<t\le1$. Thus, in order to prove \eqref{xtoinfty}, and thus the lemma, it is sufficient to show that, for $|s|=1/2$, $t\ge0$, and some positive $\delta$,
\begin{align}\label{t1p}
t+\frac{1}{\pi}\log|1-e^{2i\pi s}|>\delta.
\end{align}
>From now on we assume that $|s|=1/2$.

1) Let $0\le \arg s\le \pi/4$. Then $\sqrt{2}/2\le 2\Re s\le 1$ or $\sqrt{2}/4\le \Re s\le 1/2$. Similarly, $0\le\Im s\le \sqrt{2}/4$. Further 
\begin{align*}
-1\le\Re e^{2i\pi s}=e^{-2\pi \Im s}\cos(2\pi \Re s)\le e^{-2\pi \frac{\sqrt{2}}{4}}\cos(2\pi \frac{\sqrt{2}}{4})\le-0.026.
\end{align*}
Thus,
\begin{align*}
t+\frac{1}{\pi}\log|1-e^{2i\pi s}|\ge \frac{1}{\pi}\log|1-\Re e^{2i\pi s}|\ge\frac{1}{\pi}\log 1.026\ge0.008.
\end{align*}

2)  Let $\pi/4\le \arg s\le \pi/3$. Then  $1/4\le \Re s\le \sqrt{2}/4$, $ \sqrt{2}/4\le\Im s\le \sqrt{3}/4$, and $\Re e^{2i\pi s}\le0$. Therefore, 
\begin{align*}
t+\frac{1}{\pi}\log|1-e^{2i\pi s}|\ge t\ge \frac{\sqrt{2}}4.
\end{align*}

3)  Let $\pi/3\le \arg s\le \pi/2$. Then  $0\le \Re s\le 1/4$, $ \sqrt{3}/4\le\Im s\le 1/2$, and $\Re e^{2i\pi s}\le e^{2i\pi \sqrt{3}/4}$. This gives 
\begin{align*}
t+\frac{1}{\pi}\log|1-e^{2i\pi s}|\ge \frac{\sqrt{3}}4+\log(1-e^{-2\pi\sqrt{3}/4})\ge 0.41.
\end{align*}

4) Let $\pi/2\le \arg s\le \pi$. By a symmetry and cases 1), 2), 3) we have
\begin{align*}
t+\frac{1}{\pi}\log|1-e^{2i\pi s}|=t+\frac{1}{\pi}\log|1-e^{\overline{2i\pi s}}|\ge0.008.
\end{align*}
This proves \eqref{xtoinfty} with $\delta=0.007$ and  Lemma \ref{sn12}.

\endproof

\begin{Lemma}\label{2g-2}
Suppose $Z(s)=f(h(s))$ with an entrie function $h$ and $f(s)=a_ks^k+a_{k-1}s^{k-1}+\dots+a_1s+a_0$, where $a_k\neq 0$. Then $k$ divides $2g-2$.
\end{Lemma}
 \proof We consider the function
$$Z(s)-a_{k-1}g^{k-1}(s)-\dots-a_0=a_kg^k(s).$$
If $Z(s)$ is large, then 
$$g(s)\ll Z^\frac1k(s).$$
This, Lemma \ref{sn12}, and  Rouch\' e's theorem give  that on the circle $|s+n|=1/2$,  for large negative $n$, the functions $Z(s)$ and 
\begin{align}\label{kpower}
Z(s)-a_{k-1}g^{k-1}(s)-\dots-a_0=a_kg^k(s)
\end{align}
have the same number ($= (2g-2)(2n+1)$) of zeros (counting multiplicities). In view of formula \eqref{kpower}, for any large $n$, there is an integer $m$ such that  
$$(2g-2)(2n+1)=km.$$
Hence, $k$ divides $2g-2$.
\endproof

\proof[Proof of Theorem \ref{prime}] 
By Proposition \ref{propo} it follows that $Z$ is pseudo-prime and right-prime. 

Next we consider a decomposition $Z(s)=f(h(s))$, where $f$ is rational and $h$ is meromorphic. Since $Z$ is entire, we can assume that $h$ is entire. By Lemma \ref{pseaudoprimepolynomial} it follows from $Z(s)=f(g(s))$, with rational $f$ and entire $g$ that (i) $f$ is a polynomial or (ii) $h(z)=\exp(az+b)-w$ for some complex numbers $a\ne0$, $b$, and $w$.

In the case of (i) the theorem follows in view of Lemma \ref{2g-2}.

We shall show that case (ii) is impossible. In fact, by (ii) we have $Z(s)=f(\exp(az+b)-w)$, where $f$ is rational. It is easy to see that in the disc $|s|<R$ such a function must have less than $c_1R$ zeros (counted with multiplicities). Since  $Z(s)$ has more than $c_2 R^2$ in the disc  $|s|<R$ (as mentioned in the introduction), we arrive at a contradiction. Theorem \ref{prime} is proved. \endproof

For the proof of Corollary \ref{z1z2strong} we shall use the following Lemma due to Edrei \cite{edr}.

\begin{Lemma}\label{Edrei}
Let $f$ be an entire function. Assume that there is an unbounded sequence $(a_j)_{j=1}^\infty$ such that all but a finite number of the roots of the equations $f(z)=a_j$, $j=1,2,\dots$, lie on a straight line. Then $f$ is a polynomial of degree not greater than two.
\end{Lemma}

\proof[Proof of Corollary \ref{z1z2strong}]
We consider
\begin{align*}
\Xi_1(s)=\left((2\pi)^s\Gamma_2(s)\Gamma_2(s+1)\right)^{2g-2}Z_1(s),
\end{align*}
where the double gamma function is defined by
\begin{align*}
\frac1{\Gamma_2(s+1)}=(2\pi)^{s/2}\exp\left(-\frac{\gamma+1}{2}s^2-\frac{s}2\right)\prod_{n=1}^\infty\left(1+\frac sn\right)^n\exp\left(\frac{s^2}{2n}-s\right),
\end{align*}
and $\gamma$ is the Euler constant (see Minamide \cite{min2010}). Analogously, we define $\Xi_2(s)$. Then, for $j=1,2$, the zeros (counting  multiplicities) of $\Xi_j(s)$ coincide with the nontrivial zeros of $Z_j(s)$. By assumption, $\Xi_1(s)$ and $\Xi_2(g(s))$ have the same zeros and both are functions of order two. hence, by Hadamard's theorem, there are complex constants $a$, $b$, $c$, and $d$ such that 
\begin{align}\label{deas}
\Xi_2(g(s))=d e^{as^2+bs+c}\Xi_1(s).
\end{align}
We apply Lemma \ref{accumlines} to the function $\Xi_2(g(s))$ with $a_1=0$ and $a_2=\infty$. Clearly, $\Xi_2(g(s))$ has no poles. By \eqref{deas} its zeros lie on the line $\sigma=1/2$ apart from finitely many zeros on the real line. Thus the set $$E=\{s : \Xi_2(g(s))=a_j, j=1,2\}$$
has two accumulation lines, namely $\arg s=\pi/2$ and $3\pi/4$.  
Then Lemma \ref{accumlines} implies that  $\Xi_2(g(s))$ is a pseudo prime function. Thus $g$ must be a polynomial. By Lemma \ref {Edrei} this polynomial is of degree not greater than two. The nontrivial zeros (with only a finite number of exceptions) of both functions, $Z_1$ and $Z_2$ lie on the same line $\sigma=1/2$. Therefore $g(s)=b s$ and $b=1$. This proves Corollary \ref{z1z2strong}. \endproof

\section{Concluding Remarks}\label{cr}

If we consider a hypothetical decomposition of $Z$ with a quadratic polynomial and an entire function $g$, i.e., $Z(s)=ag(s)^2+bg(s)+c$ with complex coefficients $a\neq 0$, $b$ and $c$, then it follows that the set of $c$-points equals the union of the sets of zeros of $g$ and the $b/a$-points of $g$. More precisely, writing ${\mathcal N}(c,f)$ for the set of preimages of $f(s)=c$, we have 
$$
{\mathcal N}(c,Z)={\mathcal N}(0,g)\cup {\mathcal N}(b/a,g).
$$
Replacing $c$ by an arbitrary complex number $d$, we arrive at 
$$
{\mathcal N}(d,Z)={\mathcal N}(e,g)\cup {\mathcal N}(f,g),
$$
where $e,f$ are the solutions $(-b\pm\sqrt{b^2-4a(c-d)})/(2a)$ of the quadratic equation $aX^2+bX+c-d=0$. Taking into account the clustering of the zeros and $d$-points in general (see \cite{garunkstis12}), it appears that $g$ has to share quite a few patterns of $Z$'s value-distibution. One could expect that $g$ as well needs to be representable as a Dirichlet series in some right half-plane and has to satisfy a functional equation of Selberg-type. In view of this one could be tempted to guess that there is no non-trivial decomposition of $Z$ at all.
\bigskip

\noindent {\bf Acknowledgements.} The first author is grateful for the hospitality of Department of Mathematics of W\"urzburg University where part of this research was performed.
Moreover, the research of the first author is funded by the European Social Fund according to the activity `Improvement of researchers' qualification by implementing world-class R\&D projects' of Measure  No. 09.3.3-LMT-K-712-01-0037.

\end{document}